\documentclass{article}

\usepackage{comment}
\usepackage{fullpage}
\usepackage{extarrows}

\usepackage[usenames,dvipsnames]{pstricks}
\usepackage{epsfig}
\usepackage{pst-grad} 
\usepackage{pst-plot} 
\usepackage[space]{grffile} 
\usepackage{etoolbox} 
\makeatletter 
\patchcmd\Gread@eps{\@inputcheck#1 }{\@inputcheck"#1"\relax}{}{}
\makeatother

\usepackage{amsmath}
\usepackage{amssymb}
\usepackage{amsthm}
\usepackage{tikz,pgfplots}
\usepackage{blindtext}
\usepackage{hyperref}

\usepackage{graphicx}

\usepackage{ifthen}
\usepackage{color}

\newcommand{\intav}[1]{\mathchoice {\mathop{\vrule width 6pt height 3 pt depth  -2.5pt
\kern -8pt \intop}\nolimits_{\kern -6pt#1}} {\mathop{\vrule width
5pt height 3  pt depth -2.6pt \kern -6pt \intop}\nolimits_{#1}}
{\mathop{\vrule width 5pt height 3 pt depth -2.6pt \kern -6pt
\intop}\nolimits_{#1}} {\mathop{\vrule width 5pt height 3 pt depth
-2.6pt \kern -6pt \intop}\nolimits_{#1}}}

\def\polhk#1{\setbox0=\hbox{#1}{\ooalign{\hidewidth\lower1.5ex\hbox{`}\hidewidth\crcr\unhbox0}}}

\renewcommand{\div}{\operatorname{div}}

\newcommand{\divergence}{\operatorname{div}}

\newcommand{\dist}{\operatorname{dist}}

\newcommand{\sgn}{\operatorname{sgn}}

\newcommand{\bd}{{\bf d}}

\newtheorem{theorem}{Theorem}

\newtheorem{definition}{Definition}
\newtheorem{lemma}{Lemma}
\newtheorem{corollary}{Corollary}
\newtheorem{proposition}{Proposition}
\newtheorem{remark}{Remark}

\newcommand{\abs}[1]{\left\lvert#1\right\rvert}

\usepackage{enumerate}

\usepackage{setspace}
\begin{document}

\title{A two-phase quenching-type problem for the p-Laplacian}

\author{Julio C. Correa and Disson dos Prazeres}

\date{\today}

\maketitle

\begin{abstract}
\begin{spacing}{1.15}
\noindent We study minimizers of non-differentiable functionals of the  Alp-Phillips type with two-phases for the $p$-Laplacian , focusing on the geometric and analytical properties of free boundaries. The main result  establishes finite $(n-1)$-dimensional Hausdorff measure estimates, achieved through optimal gradient decay estimates, a $BV$-inequality and the known classifications of blow-up profiles of the linear case.
\end{spacing}

\medskip


\medskip 

\end{abstract}

\vspace{.1in}

\section{Introduction}

\subsection{The problem}
In this paper, we investigate the minimizers for functional of the type
\[
J(u)=\int_{B_1}\left(\frac{|Du|^p}p+\mathcal{F}(u)\right)\,\mathrm dx,
\]
where \( p > 2 \), and \( \mathcal{F}(u) \) is a H\"older potential defined by 
\[
\mathcal{F}(s) := \lambda_1 (s_+)^\gamma + \lambda_2 (s_-)^\gamma,
\]
where: $p>2$, $s_{\pm}=\max(\pm s,0)$, \( \lambda_1, \lambda_2 \ge 0 \) satisfying \( \lambda_1 + \lambda_2 > 0 \), and \( 0 \leq \gamma < p/2 \). Minimizers are sought in the set \( \{ u: u - g \in W^{1,p}_0(B_1) \} \), where \( g \in W^{1,p}(B_1)\cap L^\infty(B_1) \). The existence of minimizers can be established via the direct method in the calculus of variations, as the functional is coercive and weakly lower semicontinuous. However, due to the lack of convexity, the uniqueness of minimizers is not guaranteed, and multiple minimizers with the same boundary data \( g \) may exist (see \cite{Phi83}).

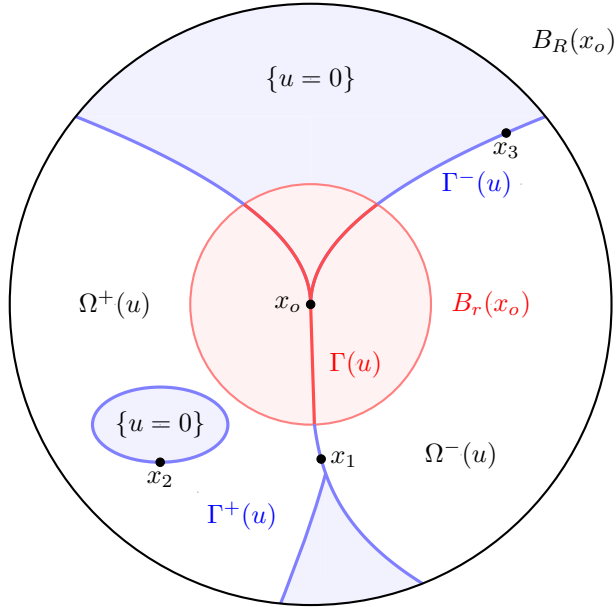
\begin{figure}[!h]\label{Figure}
\centering
\begin{tikzpicture}[scale=2]
\def\a{1/2*(1 - sqrt(17))}  
\def\b{1/2*(-1 + sqrt(17))}  
\def\c{1/2*(1 - sqrt(5))}
\def\d{1/2*(-1 + sqrt(5))}

\fill [domain=\a:0, fill=blue!5, samples=1000,  smooth]     plot (\x, {sqrt(-\x)}) |- cycle;
\fill [domain=\a:\b, fill=blue!5, smooth]     plot (\x, {sqrt(4-\x*\x)}) |- cycle;
\fill [domain=0:\b, fill=blue!5, samples=1000, smooth] plot (\x,{sqrt(\x)}) -| cycle;
\fill [fill=blue!5, samples=1000, smooth] (-0.2,-1.9899).. controls (-0.2,-1.989) and (0,-1.5) .. (0.1,-1.124) |- cycle;
\fill [domain=0.1:0.746, fill=blue!5, samples=1000, smooth] plot (\x,{-2*sqrt(sqrt(\x))}) -| cycle;
\fill [domain=-0.1:0.746, fill=blue!5, smooth]     plot (\x, -{sqrt(4-\x*\x)}) -| cycle;
domain=-0.2:0.746
\fill[fill=blue!5 , smooth] (-1,-0.8) ellipse (0.45 and 0.25);

\filldraw[color=red!50, fill=red!5, thick] (0,0) circle (0.8) ;

\draw[domain=\a:0, color=blue!50 , samples=1000,very thick] plot (\x, {sqrt(-\x)});

\draw[domain=0:\b, color=blue!50 , samples=1000,very thick] plot (\x, {sqrt(\x)});

\draw[domain=0:0.025, color=red!70 , samples=1000,very thick] plot (\x, {-2*sqrt(sqrt(\x))});

\draw[domain=0.025:0.746, color=blue!50 , samples=1000,very thick] plot (\x, {-2*sqrt(sqrt(\x))});

\draw[color=blue!50 , very thick] (-0.2,-1.9899).. controls (-0.2,-1.989) and (0,-1.5) .. (0.1,-1.124);

\draw[color=blue!50 , very thick] (-1,-0.8) ellipse (0.45 and 0.25);

\draw[domain=-0.4433:0, color=red!70 , samples=1000,very thick] plot (\x, {sqrt(-\x)});

\draw[domain=0:0.4433, color=red!70 , samples=1000,very thick] plot (\x, {sqrt(\x)});



\filldraw[black] (0,0) circle (0.8pt) node[anchor=east]{ $x_o$ };

\filldraw[black] (0.07,{-2*sqrt(sqrt(0.07))}) circle (0.8pt) node[anchor=west]{$x_1$};

\filldraw[black] (-1,-1.05) circle (0.8pt) node[anchor=north]{ $x_2$ };

\filldraw[black] (1.3,1.1401) circle (0.8pt) node[anchor=north]{ $x_3$ };

\filldraw[black] (0,1.5) circle (0.01pt) node[]{ $\{u=0\}$};

\filldraw[black] (-1.3,0) circle (0.01pt) node[]{ $\Omega^+(u)$};

\filldraw[black] (1,-1) circle (0.01pt) node[]{ $\Omega^-(u)$};

\filldraw[black] (-1,-0.8) circle (0.01pt) node[]{ $\{u=0\}$};

\filldraw[black] (-0.75,-1.25) circle (0.01pt) node[anchor=north west]{ \textcolor{blue}{$\Gamma^+(u)$}};

\filldraw[black] (1.4,0.8) circle (0.01pt) node[anchor=east]{ \textcolor{blue}{$\Gamma^-(u)$}};

\draw[thick] (0,0) circle (2);

\filldraw[black] (1.75,1.75) circle (0.01pt) node[]{$B_R(x_o)$};
\filldraw[red] (1.2,0) circle (0.01pt) node[]{$B_r(x_o)$};
\filldraw[red] (0.3,-0.4) circle (0.01pt) node[]{$\Gamma(u)$};
\end{tikzpicture}
\caption{$\Gamma^+(u):=\partial\{u>0\}\cap B_R(x_0)$, $\Gamma^-(u):\partial\{u<0\}\cap B_R(x_o)$, $x_o$ branching two-phase point, $x_1$ non-branching two-phase point, $x_2$ and $x_3$ one-phase points. The red portion represents a typical neighbourhood of a branching two-phase point.}
 
\end{figure}

A feature of our study is that no sign constraint is imposed on the boundary data \( g \). Consequently, minimizers can take both positive and negative values, leading us to regard the regions 
\[
 \Omega^+(u) := \{ u > 0 \}\hspace{.2in}\text{and}\hspace{.2in} \Omega^-(u) := \{ u < 0 \}
\]
as the distinct positive and negative phases of \( u \). Notice that since \( J \) is non-smooth,  the zero set \( \{ u = 0 \} \) may have positive Lebesgue measure. Thus, if we denote the phase interfaces as 
\[ \Gamma^\pm(u) := \partial\Omega^\pm(u) \cap B_1, \] 
 which we also call {\it free boundaries}, since they are apriori unknown; we see that $\Gamma^+(u)$ and $\Gamma^-(u)$ may split from each other. Therefore we may distinguish three type of points on these free boundaries, namely: one-phase points, where \( u \) does not change of sign sign; non-branching two-phase points, where both positive and negative phases meet but they coincide; and branching two-phase points, where the free boundaries split apart. See the Figure \ref{Figure}.

 The minimizers of the functional satisfy the Euler-Lagrange equation
\[
\Delta_p u := \mathrm{div}(|Du|^{p-2}Du) = \gamma\left( \lambda_1 (u^+)^{\gamma - 1} - \lambda_2 (u^-)^{\gamma - 1} \right)\hspace{.2in}\text{a.e in}\hspace{.1in} B_1,
\]
in the sense of the first domain variations of \( J \) (see Definition \ref{variationalsolution}). Due to the non-linear nature of the $p$-Laplacian we can not use the standard methods for the linear case.
\subsection{Known results}
The problem has been extensively studied by the profession in various settings. In the linear case (i.e. \( p = 2 \)), the problem was introduced by Alt, Caffarelli, and Friedman in \cite{AltCafFri84} for \( \gamma = 1 \), sparking significant interest in the analysis of free boundary problems of this kind. Recent works, such as \cite{PhiSpoVel18}, have culminated in the complete classification of the free boundary. For \( 0 < \gamma < 2 \), foundational results were established by Alt and Phillips in \cite{AP,Phi83}. For a comprehensive account of the one-phase setting we refer the reader to the survey \cite{Sil23}, see also \cite{SPA,T-DC}. While \cite{LinPet08} provides an extensive reference for the two-phase setting, see also \cite{PTU25}.\\

In the semi-linear case (\( p > 2 \)), the literature is less developed, primarily due to the absence of monotonicity formulae. For \( p \) close to 2 and \( \gamma = 1 \), Edquist and Lindgren in \cite{EdqLin09} proved the optimal regularity and nondegeneracy of solutions, along with estimates on the \( n-1 \)-dimensional Hausdorff measure of the free boundary near branching points. The one-phase semi-linear setting has been studied in \cite{AraTeyVos22} and \cite{AraTeyUrb24}, with results on free boundary regularity, solution nondegeneracy, and optimal growth. Some properties of the solutions of free boundary problems of obstale-type with two phases were considered for a class of heterogeneous quasilinear elliptic operators, including the $p$-Laplacian in \cite{Rod15} and more recently in \cite{CamRod24} inhomogeneous two-phase obstacle-type problem associated to the $s$-fractional $p$-Laplacian was studied. Finally, in \cite{Tei22}, Teixeira explored analogies between the loss of regularity of \( p \)-harmonic functions at critical points and the corresponding loss of regularity in free boundary problems.

\subsection{Main Results}
In this paper, we extend the classification of "blow-ups" for the linear case, as established by Weiss in \cite{Wei98}, to a more general nonlinear setting. Our main results address the optimal growth of solutions and the behavior of the free boundary in variational problems with nonlinear diffusion. The analysis of such problems is significantly more challenging due to the absence of a monotonicity formula, which plays a crucial role in the linear case. By use techniques tailored to the linear context, we provide a comprehensive understanding of the solution's growth and free boundary properties, carried in some sense by the stability of the problem.

The first main result, Theorem \ref{thm:Optimal_growth_gamma01}, establishes optimal growth estimates for solutions $u \in \mathcal{P}_1(M,p,\gamma)$ (see Definition \ref{classofminimizers})  under certain conditions on $p$ and $\gamma$ (in fact we ask for $p\sim 2$). The second main result, Theorem \ref{Thm:FinitePerimeter}, gives an $\mathcal{H}^{n-1}$ estimate on the free boundary, demonstrating that its $(n-1)$-dimensional Hausdorff measure is uniformly bounded. Our final main result, Theorem \ref{thm:Optimal_growth_dim2} concerns the two-dimensional case, in which the information on the regularity of $p$-harmonic functions together with a  Liouville-type result allows to get rid of the condition on $p$. These results are achieved by overcoming the difficulties associated with the nonlinear nature of the problem, particularly in the absence of monotonicity formulas, and provide new insights into the regularity and structure of the solutions and their free boundaries.

\begin{theorem}[Optimal Growth ]\label{thm:Optimal_growth_gamma01}
Given $0<\gamma<p/2$ and $u\in \mathcal{P}_1(M,p,\gamma,y)$. There exist positive constants $C$, $\delta$, and $r_o$, depending on $\lambda_i$, $M$, and the dimension, such that: if $\abs{p-2}<\delta$ then
\begin{equation}\label{Eq:optimal_growth}
\sup_{B_r(y)}\abs{u}\le Cr^\frac p{p-\gamma},
\end{equation}
for all $r<r_o$ and $y\in B_{1/2}\cap\Gamma(u)$.
\end{theorem}

\begin{theorem}[$\mathcal{H}^{n-1}$- Estimates on the Free Boundary]\label{Thm:FinitePerimeter}
Under the same hypotheses as Theorem \ref{thm:Optimal_growth_gamma01} and for $p>2$, we have
\begin{equation}
\mathcal{H}^{n-1}\left((\Gamma^+(u)\cup\Gamma^-(u))\cap B_{1/2}\right)\le C,
\end{equation}
where $C$ is a dimensional constant.
\end{theorem}

\begin{theorem}\label{thm:Optimal_growth_dim2}
    Given $p\ge 2$, $0<\gamma<p/2$ and $u\in \mathcal{P}_1(M,p,\gamma,y)$. There exist positive constants $C$ and $r_o$, depending on $\lambda_i$, $M$, and the dimension, such that: if $\frac{p}{p-\gamma}<1+\alpha_p$ then
\begin{equation}\label{Eq:optimal_growth n=2}
\sup_{B_r(y)}\abs{u}\le Cr^\frac p{p-\gamma},
\end{equation}
for all $r<r_o$ and $y\in B_{1/2}\cap\Gamma(u)$, where $\alpha_p$ is as in Definition \ref{def:Alphap}. 
\end{theorem}

\section{Existence of Variational Solutions}
We begin introducing the notion of a {\it variational solution} of the following equation
\begin{equation}\label{EDP}
    \Delta_pu=\gamma(\lambda_1u_+^{\gamma-1}-\lambda_2u_-^{\gamma-1})\hspace{.2in}\text{in}\hspace{.1in}B_1,
\end{equation}
where $p,\gamma,\lambda_1$ and $\lambda_2$ are constants such that: $p$ is sufficient close to $2$, $0\le \gamma\le p/2$, and $\lambda_1$  and $\lambda_2$ are non-negative constants such that $\lambda_1+\lambda_2>0$.
\begin{definition}\label{variationalsolution}
    We say that $u\in W^{1,p}_{loc}(B_1)$ is a variational solution of \eqref{EDP} provided:
    \begin{enumerate}[(i)]
    \item $u\in\mathcal{C}(B_1))\cap\mathcal{C}^{1,\alpha}((\{u>0\}\cup\{u<0\})\cap B_1)$,
    \item $(\lambda_1u_+^\gamma+\lambda_2u_-^\gamma)\in L^1_{loc}(B_1)$, 
    \item The first variation with respect to the domain variations of the functional 
    \begin{equation}\label{Functional}
    J[v]:=\int _{B_1}\left(\frac {\abs{Dv}^p}{p}+\lambda_1v_+^\gamma+\lambda_2v_-^\gamma\right)\,\bd x,
    \end{equation}
    vanishes at $u=v$.
    \end{enumerate}
\end{definition}
    \begin{remark}
        Notice that $(iii)$ above is equivalent to say that for $\xi\in\mathcal{C}^\infty_c(B_1,\mathbb{R}^n)$, we have that 
    \begin{align*}
    0=&-\left.\frac{d}{dt}\right|_{t=0}J[u(x+t\xi(x))] \\
    =&\int_{B_1}\left(\frac{|D u|^p}{p}\div\xi-p|D u|^{p-2}\left\langle D u, D\xi D u\right\rangle+(\lambda_1u_+^\gamma+\lambda_2u_-^\gamma)\div\xi\right)\,\bd x.
    \end{align*}
    \end{remark}

\bigskip
\noindent Since the functional 		
\[
I[u]:=\int_{B_1}\frac{|D u|^p}{p}\bd x,
\]
is weakly lower semicontinuous in $W^{1,p}_{loc}(B_1)$, and $(\lambda_1u_+^\gamma+\lambda_2u_-^\gamma)\in L^1_{loc}(B_1)$, we have that the dominated convergence theorem ensures that $J$ is weakly lower semicontinuous in $W^{1,p}(B_1)$ and by the condition on $\lambda_1$ and $\lambda_2$ the existence of, at least, one minimizing of \eqref{Functional} is guaranteed.\\

 Now we turn our attention to the prove of the existence of {\it variational solutions} of \eqref{EDP}.
\begin{proposition}\label{existence_of_slns}
    Let $u$ be a minimizer for $J$ then $u$ is a solution of \eqref{EDP} in the sense of the Definition \ref{variationalsolution}.
\end{proposition}

\begin{proof}
    Let $\xi\in\mathcal{C}_{c}^\infty(B_1,\mathbb{R}^n)$ be a given vector field with compact support in $B_1$ and let
    \[
    \Psi_t(x):=x+t\xi(x)\hspace{.2in} \text{for every}\hspace{.1in} x\in B_1.
    \]
    For $t$ small enough we know that $\Psi_t:B_1\to B_1$ is a diffeomorfism, so, letting $\Phi_t:=\Psi_t^{-1}$ the function 
    \[
    u_t:=u\circ\Phi_t
    \]
    is well-defined and belongs to $W^{1,p}_{loc}(B_1)$. In order to calculate the first variation with respect to the domain variations of $J$ we must calculate
    \[
     \left.\frac d{dt}\right|_{t=0}J[u_t]:=   \left.\frac d{dt}\right|_{t=0}\int_{B_1}\left(\frac{|D u_t(x)|^p}p+\lambda_1(u_t)_+^\gamma+\lambda_2(u_t)_-^\gamma\right)\, \bd x.
         \]
To do so, we will proceed into two steps.\\        

    \noindent{\bf Step 1 - }
    We begin by calculating
    \[
    \left.\frac d{dt}\right|_{t=0}\int_{B_1}\left(\lambda_1(u_t(x))_+^\gamma+\lambda_1(u_t(x))_-^\gamma\right)\ \bd x.
    \]
    Notice that, since
    \[
    x\in\{u_t>0\}\Longleftrightarrow u_t(x)>0\Longleftrightarrow \Phi_t(x)\in\{u>0\},
    \]
    we can claim that, for all $x\in B_1$
    \[
    (u_t)_+(x)=u_+(\Phi_t(x)),
    \]
    therefore, using the change of variables $y=\Phi_t(x)$, we have
    
    \[
    \int_{B_1}(u_t^+(x))^\gamma\,\bd x=\int_{B_1}(u(\Phi_t(x))_+^\gamma\, \bd x=\int_{B_1} (u(y))_+^\gamma\abs{\det \Psi_t(y)}\, \bd y.
    \]
   Since 
    \[
    |\det \Psi_t|=1+t\div \xi+o(t),
    \]
    
    we get
    \[
    \int_{B_1}(u_t)_+^\gamma\ \bd x=\int_{B_1}u_+^\gamma\ \bd x+t\int_{B_1\cap\{u>0\}}u^\gamma\div \xi\ \bd x + o(t),
    \]
    and thus 
    \[
    \left.\frac d{dt}\right|_{t=0}\int_{B_1}(u_t)_+^\gamma\,\bd x=\int_{B_1\cap\{u>0\}}u^\gamma\div\xi\ \bd x.
    \]
    
    In a similar fashion, we get
     \[
    \left.\frac d{dt}\right|_{t=0}\int_{B_1}(u_t)_-^\gamma\,\bd x=\int_{B_1\cap\{u<0\}}(-u)^\gamma\div\xi\ \bd x.
    \]
    Therefore, we have that
    \begin{equation}\label{Eq1.Prop_Existence}
    \left.\frac d{dt}\right|_{t=0}\int_{B_1}\mathcal{F}(u_t)\ \bd x=\int_{B_1}\mathcal{F}(u)\div \xi\ \bd x.
    \end{equation}
    
  \noindent{\bf Step 2 - } Now, we calculate
  \[
    \left.\frac d{dt}\right|_{t=0} \int_{B_1}\frac {|D u_t|^p}p\ \bd x. 
  \]
 By applying the change of variables $y=\Phi_t(x)$ the Inverse Function Theorem yields
 \begin{align*}
\frac 1p\int_{B_1}\left(|D u_t|^2\right)^{p/2}\ \bd x=
&\frac 1p\int_{B_1}\left(|D u(\Phi_t(x))D\Phi_t(x)|^2\right)^{p/2}\ \bd x\\
=&\frac 1p\int_{B_1}|D u(y)D\Phi_t(\Psi_t(y))|^{p/2}|\det D\Psi_t(y)|\ \bd y\\
=&\frac 1p\int_{B_1}|D u(y)[D\Psi_t(y)]^{-1}|^{p/2}|\det D\Psi_t(y)|\ \bd y\\
=&\frac 1p\int_{B_1}\left(\left\langle D u[D\Psi_t]^{-1},D u[D\Psi_t]^{-1}\right\rangle\right)^{p/2}|\det D\Psi_t|\ \bd x,
 \end{align*}
 on the other hand, since 
 \[
 D\Psi_t=Id+tD\xi,\qquad\text{and}\qquad [\Psi_t]^{-1}= Id-tD\xi+o(t),
 \]
 we have that
 \begin{align*}
 \int_{B_1}\frac{|D u_t|^p}p\ \bd x=&\frac 1p\int_{B_1}\left(\langle D u(Id-tD\xi+o(t)),D u(Id-tD\xi+o(t))\rangle\right)^\frac p2(1+t\div\xi+o(t)) \bd x\\
 =&\frac 1p\int_{B_1}\left(|D u|^2-2t\langle D u,D uD\xi\rangle\right)^\frac p2\ \bd x\\
 &+\frac 1p\int_{B_1}t\left(|D u|^2-2t\langle D u,D uD\xi\rangle\right)^\frac p2\div\xi\ \bd x+o(t), 
\end{align*}
Hence we have that
\begin{equation}\label{Eq2.Prop_Existence}
\begin{split}
\left.\frac d{dt}\right |_{t=0}\int_{B_1}\frac{|D u_t|^p}p\ \bd x=&\frac 1p\int_{B_1}\left.\frac{\partial}{\partial t}\right|_{t=0} \left(|D u|^2-2t\langle D u,D uD\xi\rangle\right)^\frac p2\ \bd x\\
&+\frac 1p\int_{B_1}\left.\frac{\partial}{\partial t}\right|_{t=0} t\left(|D u|^2-2t\langle D u,D uD\xi\rangle\right)^\frac p2\div\xi\ \bd x\\
=&-\int_{B_1}|D u|^{p-2}\langle D u,D u D\xi\rangle\ \bd x+\frac 1p\int_{B_1}|D u|^p\div\xi\ \bd x\\
=&\frac 1p\int_{B_1}|D u|^p\div\xi\ -p|D u|^{p-2}\langle D u,D u D\xi\rangle\ \bd x.
\end{split}
\end{equation}

Thus, the result follows by combining \eqref{Eq1.Prop_Existence} and \eqref{Eq2.Prop_Existence}

\end{proof}

We conclude this section with a closer look at 1-D solutions of \eqref{EDP} . To this, let  $p\ge 2$, $0<\gamma<p/2$, and consider the function 
\[
u(t):=C_1t_+^\eta-C_2t_-^\eta,
\]
where $C_i:=C_i(\lambda_i, p, \gamma)$ is a positive constant to be specified latter. Thus If $u$ where a minimizer for $J$, we must have
\begin{equation}\label{Eq1:1-D_sol}
\Delta_p u(t)=\lambda_1\gamma(C_1\eta)^{\gamma-1}(t_+)^{\eta(\gamma-1)}-\lambda_2\gamma(C_2\eta)^{\gamma-1}(t_-)^{\eta(\gamma-1)}.
\end{equation}
On the other hand, since
\[
u'(t)=C_1\eta(t_+)^{\eta-1}-C_2\eta (t_-)^{\eta-1}
\]
and
\[
\abs{u'(t)}^{p-2}u'(t)=
\begin{cases}
(C_1\eta)^{p-1}\abs{t}^{(\eta-1)(p-1)}\hspace{.4in}&\text{if}\hspace{.2in}t\ge 0\\
-(C_2\eta)^{p-1}\abs{t}^{(\eta-1)(p-1)}&\text{if}\hspace{.2in}t\le 0,
\end{cases}
\]
we have
\begin{equation}\label{Eq2:1-D_sol}
\Delta_pu(t)=(C_1\eta)^{p-1}(\eta-1)(p-1)(t_+)^{(\eta-1)(p-1)-1}-(C_2\eta)^{p-1}(t_-)^{(\eta-1)(p-1)-1}.
\end{equation}
By comparing \eqref{Eq1:1-D_sol} and \eqref{Eq2:1-D_sol}, if we let 

\begin{equation}\label{Eq3:1-D_sol}
 \eta=\frac p{p-\gamma}\hspace{.1in}\text{and}\hspace{.1in}C_i:=\frac{p-\gamma}p\left(\lambda_i\frac{p-1}{p-\gamma}\right)^{1/(p-\gamma)}\ge 2^{-1}(\lambda_i)^{1/(p-\gamma)}
\end{equation}
we can conclude that given $p,\gamma$ and $\lambda_i$ as in the hypothesis of Theorem \ref{thm:Optimal_growth_gamma01} there exists $C_i:=C(\lambda_i,p,\gamma)>0$ such that
\[
u(t):=C_1(t_+)^{p/(p-\gamma)}- C_2(t_-)^{p/p-\gamma}\hspace{.1in}\text{is a local minimizer for $J$}.
\]
Which implies the optimality in Theorem \ref{thm:Optimal_growth_gamma01}.

\section{Non-Degeneracy (Stability)}
We begin by setting some notations. Throughout the paper we will denote by $J_{B_r(y)}$ the functional
\begin{equation}\label{localfunctional}
J_{B_r(y)}[u]:=\int_{B_r(y)}\left(\frac{\abs{D u}^p}p+\lambda_1 u_+^{\gamma}+\lambda_2u_-^\gamma\right)\ \bd x.
\end{equation} 
Moreover,  we denote 
\[
\Omega^+(u):=\{u>0\}\cap B_1\hspace{.2in}\text{and}\hspace{.2in}\Omega^-(u):=\{u<0\}\cap B_1,
\]
and 
\[
\Gamma(u)=\Gamma^+(u)\cap\Gamma^-(u)\cap\{\abs{D u}=0\},
\]
where
\[
\Gamma^\pm(u):=\partial\Omega^\pm\cap B_r(y).
\]
Since our result concern only on local properties, we need to define the class of minimizers, which we are concern

\begin{definition}\label{classofminimizers}
A minimizer of the functional \eqref{Functional}  $u$ belongs to the class $\mathcal{P}_r(M,p,,\gamma,y)$ if for given $p>2$ and $0<\gamma<p/2$, the following are satisfied:
\begin{description}
\item[\text{(i)}] $u$ is a minimizer for $J_{B_r(y)}$ in the sense of Definition \ref{variationalsolution}.
\item[\text{(ii)}] $\|u\|_{L^\infty(B_r(y))}\le M$.
\item [\text{(iii)}] $y\in\Gamma(u)$
\end{description} 
\end{definition}

\begin{remark}\label{Remark:homogenityclass}
It is readily to check that under the notation above, if $u\in \mathcal{P}_{r_o}(M,p,\gamma,y)$ and $y\in\Gamma(u)$ then $u\in \mathcal{P}_r(M,p,\gamma,y)$ for $r<\dist(y,\partial B_1)$. Also we will denote $\mathcal{P}_r(M,p,\gamma)=\mathcal{P}_r(M,p,\gamma,0)$. 
\end{remark}

Our next  goal is twofold, first is to show that  the minimizers does not grow too slow around branching points, which latter will provide a non-degerency of the solutions. The second part will be  concerned to the stability of the class $\mathcal{P}_r(M,p,\gamma,y)$ as $p$ tends to 2.  As usual, constants standing for C may change from line to line, and depend only on the suitable universal quantities as $\lambda_i$, $M$ and the dimension.  The first part of the aforementioned is the following

\begin{proposition}\label{Prop:nondegeneracygamma01}
Let $p\ge2$, $0<\gamma<p/2$, $0<r_o<1$ given and $u\in \mathcal{P}_{r_o}(M,p,\gamma)$. Then for every $y\in\Gamma(u)$  given, there exists a universal constant $C>0$ such that
\begin{equation}\label{Eq1.Prop_nondegeneracygamma01}
\sup_{B_r(y)\cap\Omega^+(u)}u\ge Cr^{\frac p{p-\gamma}},
\end{equation}
and 
\begin{equation}\label{Eq2.Prop_nondegeneracygamma01}
\inf_{B_r(y)\cap\Omega^-(u)} u\le Cr^{\frac p{p-\gamma}},
\end{equation}
for all $r<\dist(y,\partial B_{r_o})$.
\end{proposition}
\begin{proof}
First, notice that once we have prove \eqref{Eq1.Prop_nondegeneracygamma01} it is rather easy to prove \eqref{Eq2.Prop_nondegeneracygamma01}, in a very similar fashion in fact, furthermore by continuity it is enough to prove the result for $y\in\Omega^{\pm}(u)\cap B_{1}$. Let $C>0$ and $\eta>0$ be a fixed constants, to be specified latter, and define the function
\[
w(x):=C\abs{x-y}^{\frac p{p-1}}\hspace{.2in}\text{in}\hspace{.1in}\Omega^{+}(u)\cap B_{1}.
\]
A direct calculation shows that, for $x\in\Omega^+(u)\cap B_{r_o}(y)$
\begin{equation}\label{Eq1:Prop_nondegeneracy01}
\begin{split}
\Delta_p(u^\eta(x))
=&\div(\abs{D(u^\eta(x))}^{p-2}D(u^{\eta}(x)))\\
=&\eta^{p-1}\div(u^{(\eta-1)(p-1)}\abs{Du(x)}^{p-2}Du(x))\\
=&\eta^{p-1}(\eta-1)(p-1)u^{(\eta-1)(p-1)-1}(x)\abs{Du(x)}^{p-2}\langle Du(x),Du(x)\rangle\\
&+\eta^{p-1}u^{(\eta-1)(p-1)}(x)\div(\abs{Du(x)}^{p-2}Du(x))\\
&=\eta^{p-1}u^{(\eta-1)(p-1)-1}(x)\left((\eta-1)(p-1)\abs{Du(x)}^p+u(x)\Delta_pu(x)\right),
\end{split}
\end{equation}
and 
\begin{equation}\label{Eq2:Prop_nondegeneracy01}
\Delta_pw(x)= C^{p-1}n\left(\frac p{p-1}\right)^p.
\end{equation}
Thus, if we choose 
\[
\eta=\frac{p+\gamma}{p-1}\in [1,3]\hspace{.2in}\text{provided}\hspace{.1in}0<\gamma<p/2,
\]
Proposition \ref{existence_of_slns}, together with the fact that $x\in\Omega^+(u)\cap B_{r_o}(y)$ and \eqref{Eq1:Prop_nondegeneracy01} yields
\begin{align*}
\Delta_p (u^\eta)
&=\eta^{p-1}u{-\gamma}\left((1-\gamma)\abs{Du}^p+u\Delta u\right)\\
&=\eta^{p-1}u^{-\gamma}\left((1-\gamma)\abs{Du}^p+u\gamma\lambda_1u^{\gamma-1}\right)\\
&=\eta^{p-1}u^{-\gamma}\left((1-\gamma)\abs{Du}^p+\lambda_1\gamma u^{\gamma}\right)\\
&=\eta^{p-1}(1-\gamma)\frac{\abs{Du}^p}{u^{\gamma}}+\lambda_1\gamma\eta^{p-1}\\
&\ge\lambda_1\gamma\eta^{p-1}>0.
\end{align*}
On the other hand, if we choose $C>0$ such that
\begin{align*}
C^{p-1}n\left(\frac p{p-1}\right)^p&\le\lambda_1\gamma\eta^{p-1}\\
C^{p-1}&\le\left(\frac{\lambda_1\gamma(p-1)}{np}\right)\eta^{p-1}\left(\frac{p-1}{p}\right)^{p-1}\\
&\le\left(\frac{\lambda_1\gamma(p-1)}{np}\right)\left(\frac{p-\gamma}{p}\right)^{p-1}\\
C&\le\left(\frac{\lambda_1\gamma(p-1)}{np}\right)^{\frac 1{p-1}}\frac{p-\gamma}p,
\end{align*}
from \eqref{Eq1:Prop_nondegeneracy01} we have
\[
\Delta_p (u^\eta) \geq \Delta_p w \hspace{.2in}\text{in}\hspace{.1in}\Omega^{+}(u) \cap B_r(y).
\]
 Since $u^\eta(y)>w(y)=0$ there is, by the comparison principle, $x_y \in\partial\left(\Omega^{+}(u)\cap B_r(y)\right)$ such that
 \[
  u^\eta\left(x_y\right)>w\left(x_y\right).
  \]
   Moreover, $u^\gamma \leq w$ on $\Gamma^{+}(u)$ so $x_y \in \partial B_r \cap \Omega^{+}(u)$. Therefore,
$$
\sup _{\partial B_r(y) \cap \Omega^{+}(u)} u>Cr^{\frac p{p-\gamma}}$$
\end{proof}
In order to tackle the stability matter, the following lemma is needed.

\begin{lemma}\label{lemma_convergence_sets}
Let $0<\gamma<p/2$ be fixed and $u_k\in \mathcal{P}_r(M,p_k,\gamma)$, where $p_k\to2$, be such that $u_k\to u$ in $\mathcal{C}^{1,\alpha}(B_{r_o})$ for some $\alpha>0$. Then
\begin{description}
\item[\text{(i)}] $\limsup_{k\in\mathbb{N}}\{u_k=\abs{u_k}=0\}\subset\{u=\abs{u}=0\}$,
\item[\text{(ii)}]$\Omega^\pm(u)\subset\liminf_{k\in\mathbb{N}}\Omega^\pm(u_k)$,
\item[\text{(iii)}]$\limsup_{k\in\mathbb{N}}\Gamma^\pm(u)\subset\Gamma^\pm(u)$.
\end{description}
\end{lemma}

\begin{proof}
Notice that {\bf (i)}  and {\bf (ii)} follows from the $\mathcal{C}^1$ convergence. In order to prove {\bf (iii)} let $(x_k)_{k\in\mathbb{N}}\subset\Gamma(u_k)$ and consider a convergent subsequence, by abusing of notation also labeled with $x_k$, such that $x_k\to x$. First, notice that by continuity $u_k(x_k)\to u(x)$  and therefore $x\in\Gamma (u)$.\\

On the other hand, the Propositions  \ref{Prop:nondegeneracygamma01}   implies that for any $0<r<\dist(x_k,\partial B_{r_o}(x))$ there exists an $y_k\in\partial B_r(x_k)$ such that
\[
u_k(y_k)\ge Cr^{\frac p{p-\gamma}}>0\hspace{.2in}\text{for all }k,
\]
where $C>0$ is an universal constant (moreover, $C$ does not depend on $k$). Thus, by continuity, there exists $y\in\partial B_r(x)$ such that
\[
u(y)\ge Cr^{\frac p{p-\gamma}}>0.
\]
 
Since $r$ is arbitrary, we can conclude that $x\in\Gamma^+(u)$, and consequently
\[
\limsup_{i\in\mathbb{N}}\Gamma^+(u_k)\subset\Gamma^+(u).
\]
In very same fashion can be shown the other inclusion.
\end{proof}
We close this section with the an stability result, which reads
\begin{proposition}\label{stability}
Let $(u_k)_{k\in\mathbb{N}}$ be a sequence such that for every $k\in\mathbb{N}$, $u_k \in \mathcal{P}_{r_o}(M,p_k,\gamma)$ and $p_k\to2$. Then there exists a subsequence, without loss of generality labeled again $u_k$ and a function $u_2\in \mathcal{P}_r(M,2,\gamma)$ such that 
\[
u_k\to u_2\hspace{.2in}\text{in}\hspace{.2in}\mathcal{C}^{1,\alpha}(B_{r_o}),
\] 
for all $r<r_o$.
\end{proposition}
\begin{proof}
If $0<\gamma<1$ then by the Theorem 1.1 in \cite{LQT} $u_k\in\mathcal{C}^{1,\alpha}(B_r)$ uniformly. To $1\leq \gamma<p-2$, since by Proposition \ref{existence_of_slns} we have that
\[
|\Delta_{p_k}u_k|\le C,
\]
in distribution sense, combined with in \cite[Theorem 1]{Lie88} yields,
\[
u_k\in\mathcal{C}^{1,\alpha}(B_r)\hspace{.2in}\text{uniformly}.
\]
Finally, the uniform convergence above, together with Lemma \ref{lemma_convergence_sets} implies that, up to a subsequence, $(u_k)_{k\in\mathbb{N}}$ converges in $\mathcal{C}^{1,\alpha}(B_r)$ to some $u_2\in\mathcal{P}_1(M,2,\gamma)$.
\end{proof}

\section{Optimal growth and its consequences}

The proof of Theorem \ref{thm:Optimal_growth_gamma01} relies on two key ingredients: First, the optimal regularity of $u_0\in\mathcal{P}_r(M,2,\gamma)$, which due to the  work of Giaquinta and Giusti  in \cite{GiaGiu84} we known that $u_0$ is $C^{\frac{2}{2-\gamma}}$ (in fact, in \cite{GiaGiu84} this regularity was proved for a general class of functionals that includes $J$. For a more direct proof of this, specify for the minimizers of $J$ see the note of Lindgren and Silvestre \cite{LinSil05}). The second ingredient is the Blow-up classifications of minimizers in the class $\mathcal{P}_r(M,2,\gamma)$ due to Weiss in \cite{Wei98}. 

\begin{proof}[Proof of Theorem \ref{thm:Optimal_growth_gamma01}]
    Without loss of generality, we will consider $y=0$. Given sequences $u_k\in \mathcal{P}_1(M,p_k,\gamma)$, $\delta_k\rightarrow 0$, we have $p_k\rightarrow 2$, and by Lemma \ref{lemma_convergence_sets}, we have $u_k\rightarrow u_0\in \mathcal{P}_r(M,2,\gamma)$ in $C^{1,\alpha}$ for all $r<1$, where $u_0$ satisfies
    $$
   \Delta u_0=\gamma\left(\lambda_1 (u_0^+)^{\gamma-1}-\lambda_2(u_0^-)^{\gamma-1}\right)\hspace{.2in}\text{in}\hspace{.1in}B_1.
    $$
   Therefore, for $0<\gamma<1$ we have,
    \begin{equation*}\label{estimates c2}
       \sup_{B_r}|u_0|\leq Cr^{\frac{2}{2-\gamma}}
    \end{equation*}
    and, for $1\le\gamma<p/2$ 
    \begin{equation*}
    \sup_{B_r}|u_0|\leq Cr^2,
    \end{equation*}
    for $r \leq 1/2$. In any case, for $k$ large enough, we have that
    $$
    \sup_{B_{r}}|u_k| \leq  \sup_{B_r}|u_k-u_0|+ \sup_{B_r}|u_0|,
    $$
    since that for $k$ large $\sup_{B_r}|u_k-u_0|<\epsilon$ we choose $r=\bar r$ in such way that
    $$
    \sup_{B_{\bar r}}|u_k|\leq C\bar r^{\frac{p_k}{p_k-\gamma}}, 
    $$
    which is possible since $p_k>2$. We consider the set
    $$
    A_k=\{ r>0; \sup_{B_{r}}|u_k|> (C+1)r^{\frac{p_k}{p_k-\gamma}} \}
    $$
    and observe that the set $A_k$ is bounded. We then consider $r_k=\sup A_k$. If we can prove that $r_k=0$ for sufficiently large $k$, then the theorem is proved. Suppose this is not true; then there exists a subsequence $r_k>0$ such that $r_k\rightarrow r_0$. We have $r_0=0$, in fact, since
    $$
    \sup_{B_{r_k}}|u_k|=(C+1)r_k^{\frac{p_k}{p_k-\gamma}}
    $$
    if $r_0>0$, by \eqref{estimates c2}, we have
    $$
     \sup_{B_{r_0}}|u_0|=(C+1)r_0^{\frac{2}{2-\gamma}}>Cr_0^{\frac{2}{2-\gamma}}\geq  \sup_{B_{r_0}}|u_0|,
    $$
    which is a contradiction.
    Now we define
    $$  
    v_k(x)=\frac{u_k(r_kx)}{r_k^{\frac{p_k}{p_k-\gamma}}}\quad \mbox{ in   } B_{1/r_k},
    $$
    thus $\sup_{B_{1}}|v_k|=C+1$ and $\sup_{B_{r}}|v_k|\leq (C+1)r^{\frac{p_k}{p_k-\gamma}}$ for $1<r<\frac{1}{r_k}$, since $r_k=\sup A_k$. Furthermore, $v_k\in P_r((C+1)r^{\frac{p_k}{p_k-\gamma}},p_k,\gamma)$ and therefore $v_k\rightarrow v_0$ locally in $C^{1,\alpha}(\mathbb{R}^n)$. Therefore $ \sup_{B_{1}}|v_0|=(C+1)$ , $\sup_{B_{r}}|v_0|\leq (C+1)r^{\frac{2}{2-\gamma}}$ for $1<r$ and $v_0\in P_r((C+1)r^{\frac{2}{2-\gamma}},2)$ for $r\geq 1.$ Therefore, $v_0$ satisfies the conditions of Theorem 4.1 in \cite{Wei98}, thus $v_0(rx)=r^{\frac2{2-\gamma}}v_0(x)$. Hence, since
     $$
     r_k^{\frac{p_k}{p_k-\gamma}}|v_k(r_k^{-1}x)|=|u_k(x)|\leq M,
     $$
     passing to the limit we have that
     $$
     \sup_{B_1}|v_0(x)|\leq M,
     $$
     which is a contradiction if we take $C>M$.
        
\end{proof}

\begin{corollary}\label{Cor:growth_gradient}
Let $u\in \mathcal{P}_1(M,p,\gamma)$, where $\abs{p-2}<\delta$ and $\delta$ as in the Theorem \ref{thm:Optimal_growth_gamma01}. There exists positive constants $C$ and $r_o$ such that
\begin{equation}\label{Eq0_cor:growth_gradient}
\sup_{B_r}\abs{D u}\le Cr^{\frac\gamma{p-\gamma}},
\end{equation}
for all $r<r_o$
\end{corollary}
\begin{proof}
    Assume the assertion fails. Then for a fixed positive constant $C_o$, to be chosen latter, we have that for all $k\in\mathbb{N}$
    \begin{equation}\label{Eq1_cor:growth_gradient}
        \mu_{k+1}\ge\max(C_o(2^{-(k+1)})^{\frac\gamma{p-\gamma}},2^{-\frac\gamma{p-\gamma}}\mu_k),
    \end{equation}
    where 
    \[
        \mu_k=\sup_{B_{2^{-k}}}\abs{D u}.
    \]

\noindent    Now, letting
        \[
        v(x):=\frac{u(2^{-k}x)}{2^{-k}\mu_{k+1}}\hspace{.2in}\text{for}\hspace{.1in} x\in B_1,
    \]
   the  Theorem \ref{thm:Optimal_growth_gamma01} and yields the existence of a universal constant $C>0$ such that
   \[
   \sup_{x\in B_1}\abs{u(2^{-k}x)}=\sup_{B_{2^{-k}}}\abs{u(x)}\le C(2^{-k})^\frac p{p-\gamma},
   \]
   and by \eqref{Eq1_cor:growth_gradient}, for $x\in B_1$ we have
   \begin{equation}\label{Eq2_cor:growth_gradient}
   \begin{split}
   \abs{v(x)}&\le \frac C{\mu_{k+1}}\frac{(2^{-k})^{\frac p{p-\gamma}}}{2^{-k}}\\
   &= \frac C{\mu_{k+1}}(2^{-k})^{\frac\gamma{p-\gamma}}\\
   &\le \frac C{C_o}2^{\frac\gamma{p-\gamma}},
   \end{split}
   \end{equation}
   thus, by choosing $C_o\ge Ck2^{\gamma/p-\gamma}$, we get
   \[
   \abs{v(x)}\le \frac 1k\hspace{.2in}\text{for}\hspace{.1in}x\in B_1.
   \] 
   Moreover
   \begin{align*}
   \sup_{B_{1/2}}\abs{D v}&=\sup_{B_{1/2}}\frac{\abs{D (u(2^{-k}x))}}{2^{-k}\mu_{k+1}}\\
   &=\mu_{k+1}^{-1}\sup_{B_{1/2}}\abs{D u(2^{-k}x)}\\
   &=1
   \end{align*}
   and, again by \eqref{Eq1_cor:growth_gradient},
   \begin{align*}
   \sup_{B _1}\abs{D v}&=\sup_{B_1}\frac{\abs{D (u(2^{-k}x))}}{2^{-k}\mu_{k+1}}\\
   &=\frac{\mu_k}{\mu_{k+1}}\\
  & \le \frac {\mu_k}{2^{-\frac\gamma{p-\gamma}}\mu_k}\\
  &=2^{\frac\gamma{\gamma-p}}.
   \end{align*}
To reach a contradiction, notice that 
\begin{align*}
\Delta_pv(x)&=\mu_{k+1}^{1-p}2^{-k}\Delta_p u(2^{-k}x)\\
&\le\mu_{k+1}^{1-p}2^{-k}C_1\abs{u(2^{-k}x)}^{\gamma-1}\\
&\le C_1\mu_{k+1}^{\gamma-p}(2^{-k})^{\gamma-1}\abs{v(x)}^{\gamma-1},
\end{align*}
where $C_1>0$ is an universal constant depending only on $\lambda_1$, $\lambda_2$ and $\gamma$. Therefore, \eqref{Eq1_cor:growth_gradient} and the choice on $C_0$  yields
\begin{align*}
    \Delta_p v&\le 2C_1C^{\gamma-p}k^{\gamma-p}\abs{v}^{\gamma-1}\\
    &\le 2C_1C^{\gamma-p}k^{1-p}.
\end{align*}
Thus, the $C^1$-estimates for the $p$-Laplacian (see \cite{Tol84}) implies that 
\[
    1=\sup_{B_1/2}\abs{D v}\le C_2\left(\sup_{B_1}v+k^{1-p}\right)\le C_2(k^{-1}+k^{1-p}),
\]
which for $k$ large enough gives a contradiction.
\end{proof}
\begin{lemma}\label{lemma_l2estimate_hessian}
    Assume that the hypothesis of  Theorem \ref{thm:Optimal_growth_gamma01}  are in force. Then there exist $C$ and $r_o$ positive, universal constants such that
    \begin{equation}\label{Eq0:lemma_l2estimate_hessian}
        \frac{1}{\abs{B_r}}\int_{B_r}\left(\abs{D u}^{p-2}\abs{D^2u}\right)^2\bd x\le C,
    \end{equation}
    for all $r<r_o$.
\end{lemma}
\begin{proof}
    Let us first notice that if we let
    \[
        \mathcal{S}_k(u)=\int_{B_1}\left(\abs{D u(2^{-k}x)}^{p-2}\abs{D^2u(2^{-k}x)}\right)^2\bd x,
    \]
    it suffices to prove that there exists a positive universal constant $C_o$ such that for all $k\in\mathbb{N}$, the following holds
    \[
        \mathcal{S}_{k+1}(u)\le\max (C_0,\mathcal{S}_{k}(u)).
    \]
    Suppose, contrary to our claim, that for given $k\in\mathbb{N}$,  there exists an integer $j_k$ such that
    \begin{equation}\label{Eq1:lemma_l2estimate_hessian}
        \mathcal{S}_{j_k+1}(u)>\max(4^{p(j_k+1)}k^2 ,\mathcal{S}_{j_k}(u)).
    \end{equation}
    Now, for such $k$, define  
    \[
    v_k(x)=\left(2^{-    p(j_k+1)}\sqrt{\mathcal{S}_{j_k+1}(u)}\right)^{-\frac{1}{p-1}}u(2^{-j_k}x)\hspace{.2in}\text{for}\hspace{.1in}x\in B_1,
    \]
    thus, from \eqref{Eq1:lemma_l2estimate_hessian} and the fact that $u\in\mathcal{P}_1(M,p,\gamma)$, it readily follows that:\\
    
    \begin{enumerate}[(i)]
        \item $\|v_k\|_{L^\infty(B_1)}\le Mk^{-\frac 1{p-1}}$\\
       \item $\mathcal S_1(v_k)=\left[\left(2^{-p(j_k+1)}\sqrt{\mathcal S_{j_k+1}(u)}\right)^{\frac{p-1}{1-p}}2^{-p(j_k+1)}\right]^2\mathcal S_{j_k+1}(u)=1$
        \item $\Delta_p v_k(x)=\left(2^{-(j_k+p)}\sqrt{\mathcal{S}_{j_k+1}(u)}\right)^{-1}\Delta_p u(2^{-j_k}x)$.
    \end{enumerate}
    
    To get  a contradiction we proceed in two steps. We begin by use the fact that for given $k\in\mathbb{N}$, $v_k\in\mathcal{P}(Mk^{-\frac1{p-1}},p,\gamma)$.\\
    
    \noindent{\bf Step - 1.} First of all, notice that from (i) and (iii) above, letting $p^*=p/(p-1)$, we get
    \begin{align*}
        \Delta_p v_k(x)&\le C_1(\lambda_i,\gamma)2^{(j_k+p)}\left(\sqrt{\mathcal{S}_{j_k+1}(u)}\right)^{-1}\abs{u(2^{-j_k}x}^{\gamma-1}\\
        &\le C_2(\lambda_1,\gamma,p)2^{(j_k+1)}2^{-p^*(\gamma-1)(j_k+1)}\left(\sqrt{\mathcal{S}_{j_k+1}(u)}\right)^{\frac{\gamma-1}{p-1}-1}\abs{v_k(x)}^{\gamma-1}\\
        &\le C_2(\lambda_1,\gamma,p)2^{(j_k+1)}2^{-p^*(\gamma-1)(j_k+1)}\left(\sqrt{\mathcal{S}_{j_k+1}(u)}\right)^{\frac{\gamma-1}{p-1}-1}2^{-p^*(p-\gamma)(j_k+1)}k^{-\frac{p-\gamma}{p-1}}M^{\gamma-1}\\
        &\le C_3(\lambda_i,\gamma,p,M)(2^{j_k})^{1-p}k^{-\frac{p-\gamma}{p-1}},
    \end{align*}
    hence, there exists an universal constant $C$ such that 
    \begin{equation}\label{Eq2:lemma_l2estimate_hessian}
       \|\Delta_p v_k\|_{L^{\infty}(B_1)}\le C\left(2^{j_k}\right)^{1-p}k^{-\frac{p-\gamma}{p-1}}.
    \end{equation}
    Furthermore, since $p>2$ the uniform $\mathcal{C}^1$-estimates (see \cite{Tol84}) yields
    \begin{equation}\label{Eq3:lemma_l2estimate_hessian}
    \begin{split}
        \sup_{B_{1/2}}\abs{D v_k}&\le C(\|v_k\|_{L^\infty(B_1)}+\|\Delta_p v_k\|_{L^\infty(B_1)})\\
        &\le C\left(k^{-\frac1{p-1}}+\left(2^{j_k}\right)^{1-p}k^{-\frac{p-\gamma}{p-1}}\right)
    \end{split}
    \end{equation}

    \noindent{\bf Step - 2.} Finally, since $p>2$, the $L^2$ bounds of the second derivatives (see \cite{Tol84}) yields, for all $k\in\mathbb{N}$,
    \begin{equation}
    \begin{split}
    S_1(v_k)&=\int_{B_{1/2}}\left( \abs{D v_k}^{p-2}\abs{D^2 v_k}\right)^2\bd x\\
    &\le\|D v_k\|_{L^\infty(B_{1/2})}^{2(p-2)}\|D^2v_k\|^2_{L^2(B_{1/2})}\\
    &\le C\left(k^{-\frac1{p-1}}+(2^{j_k})^{1-p}k^{-\frac{p-\gamma}{p-1}}\right)^{2(p-2)}\\
    &<1,
    \end{split}
    \end{equation}
    which contradicts (ii), and then the proof is complete.
\end{proof}

\section{$\mathcal{H}^{n-1}$-estimates on the free boundary}
In the sequel, we examine the local properties of  $\Gamma^+(u)\cup\Gamma^-(u)$ and present the proof of Theorem \ref{Thm:FinitePerimeter}. The following lemmas will provide an adequate estimate for the measure of the level sets (which will play a major role later) by performing integration by parts with a convenient choice of test functions. To simplify the notation, from now on, let
\[
\mathcal{F}(u):=\gamma\left(\lambda_1(u^+)^{\gamma-1}-\lambda_2(u^-)^{\gamma-1}\right),
\]
and
\[
\psi_\varepsilon:=
\begin{cases}
    1\hspace{.2in}&;\, t>\varepsilon^{\frac 1{p-\gamma}}\\
    \frac{\abs{t}^{p-\gamma}\sgn t}{\varepsilon}&\,\abs{t}\le\varepsilon^{\frac 1{p-\gamma}}\\
    -1\hspace{.2in}&;\, t<-\varepsilon^{\frac 1{p-\gamma}}\\
\end{cases}
\]
\begin{lemma}\label{lemabvinequality}
Assume that the hypothesis of Theorem \ref{Thm:FinitePerimeter} are in force, and take $r_o$ as in Corollary \ref{Cor:growth_gradient} and Lemma \ref{lemma_l2estimate_hessian}. Then for any non-negative $\varphi\in\mathcal{C}^\infty_0(B_r)$ and any $r<r_o$ we have 
\begin{equation}\label{Eq:lemmabvinequality}
    \int_{B_r}\abs{D\mathcal{F}(u)}\varphi\,\bd x=\int_{B_r}\abs{D\Delta_pu}\varphi\,\bd x\le C\int_{B_r}\abs{D\varphi}\bd x.
\end{equation}
Here $C>0$ is an universal constant, and the inequality has to be understood in the BV-sense. Besides, for given $\varepsilon>0$, there exists a universal constant C (which apriori is different from the one above) such that
    \begin{equation}\label{Eq0:lemmabvinequality}
    \begin{split}
        \int_{B_r}&\varphi\psi'_\varepsilon(\partial_i u)\left(\abs{D u}^{p-2}\abs{D\partial_i u}^2+(p-2)\abs{D u}^{p-4}\langle\partial u,D u\rangle^2\right) \bd x\\
       & \le C\int_{B_r}\abs{D\varphi}\,\bd x
    \end{split}
    \end{equation}
\end{lemma}

\begin{proof}
    We split the argument into three steps. Since the inequality in \eqref{Eq:lemmabvinequality} has to be proved in BV-sense, we start by cook-up our test functions.\\

    \noindent{\bf Step 1-}
    Let $\mathcal{F}_\varepsilon(u)$ be a non-decreasing smooth aproximation of $\mathcal{F}(u)$  and for given $0<r<r_o$ (with $r_o$ as in the hypothesis) let  $u_\varepsilon$ to be the solution of
    \[
    \begin{cases}
        \Delta_p u_\varepsilon=\mathcal{F}_\varepsilon(u)\hspace{.2in}&\text{in}\hspace{.1in} B_r\\
        u_\varepsilon=u &\text{on}\hspace{.1in}\partial B_r.
    \end{cases}
    \]
    Since 
    \[
  D(\Delta_p u_\varepsilon)=D\mathcal{F}_\varepsilon(u),
    \]
    multiplying by $\psi_\varepsilon(\partial_iu_\varepsilon)\varphi$ and integrating by parts, we get
    \begin{align*}
        \int_{B_r}\psi_\varepsilon(\partial_iu_\varepsilon)\varphi D\eta_\varepsilon(u) \bd x=&\int_{B_r}\psi_\varepsilon(\partial_iu_\varepsilon)\varphi D(\Delta_p u_\varepsilon)\bd x\\
        =&-\int_{B_r}\Delta_pu_\varepsilon D(\psi_\varepsilon(\partial_iu_\varepsilon)\varphi) \bd x\\
        =&-\int_{B_r}\Delta_pu_\varepsilon(\psi'_\varepsilon(\partial_iu_\varepsilon)D\partial_iu_\varepsilon+\psi_\varepsilon(\partial_iu_\varepsilon)D\varphi) \bd x\\
        &=-\int_{B_r}\Delta_pu_\varepsilon\psi'_\varepsilon(\partial_iu_\varepsilon)\varphi D\partial_iu_\varepsilon \bd x-\int_{B_r}D\varphi\Delta_pu_\varepsilon\psi(\partial_iu_\varepsilon) \bd x\\
        &= I-II,
    \end{align*}
    where 
    \[
    I=-\int_{B_r}\Delta_pu_\varepsilon\psi'_\varepsilon(\partial_iu_\varepsilon)\varphi D\partial_iu_\varepsilon \bd x,
    \]
    and
    \[
    II=\int_{B_r}D\varphi\Delta_pu_\varepsilon\psi(\partial_iu_\varepsilon) \bd x.
    \]
\noindent
    Furthermore, since $\psi_\varepsilon$ is non-decreasing, we get
    \[
    \abs{\int_{B_r}\psi_\varepsilon(\partial_iu_\varepsilon)\varphi D\eta_\varepsilon(u) \bd x}\le\abs{II},
    \]
    and subsequently, 
    \begin{align*}
        \abs{\int_{B_r}\psi_\varepsilon(\partial_iu_\varepsilon)\varphi D\eta_\varepsilon(u) \bd x}\le&\int_{B_r}\abs{D\varphi}\abs{\Delta_pu_\varepsilon}\abs{\psi_\varepsilon(\partial_iu_\varepsilon)} \bd x\\
        \le&\int_{B_r}\abs{D\varphi}\abs{\divergence(\abs{D u_\varepsilon}^{p-2}D u_\varepsilon)} \bd x\\
        \le&\int_{B_r}\abs{D\varphi}\sum_i\abs{\partial_i(\abs{D u_\varepsilon}^{p-2}D u_\varepsilon)} \bd x\\
        =&\int_{B_r}\abs{D\varphi}\sum_i\left(\abs{\partial_i\abs{D u}{p-2}}\abs{D u_\varepsilon}+\abs{D u_\varepsilon}^{p-2}\abs{\partial_i\abs{D u_\varepsilon}}\right)\ \bd x\\
        =&\int_{B_r}\abs{D\varphi}\sum_{i}\left((p-2)\abs{D u_\varepsilon}^{p-3}\abs{\partial_i\abs{D u_\varepsilon}}\abs{D u_\varepsilon}\right)\ bd x\\
        =&\int_{B_r}\abs{D\varphi}\left((p-1)\abs{u_\varepsilon}^{p-2}\sum_i\abs{\partial_i\abs{D u_\varepsilon}}\right)\\
        \le& C\int_{B_r}\abs{D\varphi}\abs{D u_\varepsilon}^{p-2}\abs{D^2u_\varepsilon}\ \bd x,
    \end{align*}
    which together with \eqref{Eq0:lemma_l2estimate_hessian}, yields 
    \begin{equation}
        \abs{\int_{B_r}\psi_\varepsilon(\partial_iu_\varepsilon)\varphi D\eta_\varepsilon(u) \bd x}
        \le C\int_{B_r}\abs{D\varphi}.
    \end{equation}

    \noindent{\bf Step 2 -} Now we observe that letting $\varepsilon\to 0$,  $\psi_\varepsilon$ converges to the $\sgn$ function, which combined with the $W^{2,2}$-estimates of the $p$-Laplacian ( see \cite{Tol84}) yields a subsequence, by abusing the notation labeled again as, $u_\varepsilon$ such that
    \[
        u_\varepsilon\rightharpoonup u\hspace{.2in}\text{weak in}\hspace{.1in}W^{2,2}(B_r),
    \]
    and
    \[
        u_\varepsilon\to u\hspace{.2in}\text{strongly in}\hspace{.1in} \mathcal{C}^{1,\alpha}(B_r).
    \]
    On the other hand, since $\mathcal{F}_\varepsilon(u)\to\mathcal{F}(u)$ in $H^1_0(B_r)$, the upper semi-continuity of the BV norm provide,
    \begin{align*}
        \int_{B_r}\abs{D\eta(u)}\varphi\bd x
        =&\int_{B_r}\abs{\int_{B_r}(\sgn (\partial_iu)D\eta(u))\varphi\ \bd x}\\
        \le&\lim_{\varepsilon\to0}\abs{\int_{B_r}\psi_\varepsilon(\partial_i u_\varepsilon)D\eta_\varepsilon(u)\varphi\ \bd x} \\
        \le& C\int_{B_r}\abs{D \varphi}\ \bd x,    
    \end{align*}
    which conclude the proof of \eqref{Eq:lemmabvinequality}.\\
    
    \noindent{\bf Step 3 -} Finally, we take $\psi_\varepsilon'(\partial_iu)\varphi$ as a test function, to get
    \begin{align*}
    \int_{B_r}&\left(\partial_i\Delta_p u\right)\psi'_\varepsilon(\partial_iu)\varphi\ \bd x\\
    =&-\int_{B_r}\partial_i\left(\abs{D u}^{p-2}\abs{D u}\right)D\left(\psi_\varepsilon(\partial_i u)\varphi\right)\ \bd x\\
    =&-\int_{B_r}\left(\abs{D u}^{p-2}D\partial_iu+(p-2)\abs{D u}^{p-4}\langle D\partial_i u,D u\rangle\right)\psi_\varepsilon(\partial_i u)D\varphi\,\bd x\\
    &-\int_{B_r}\left(\abs{D u}^{p-2}D\partial_iu+(p-2)\abs{D u}^{p-4}\langle D\partial_i u,D u\rangle\right)\psi'_\varepsilon(\partial_i u)\varphi D\partial_i u\,\bd x.
    \end{align*}
    Therefore, 
    \begin{align*}
        \int_{B_r}&\varphi\psi_\varepsilon'(\partial_i u)\left(\abs{D u}^{p-2}\abs{D\partial_i u}^2+(p-2)\abs{D u}^{p-4}\langle D\partial_i u,D u\rangle^2\right)\,\bd x\\
        =&-\int_{B_r}\left(\partial_i\Delta_pu\right)\psi_\varepsilon(\partial_i u)\varphi\,\bd x\\
        &-\int_{B_r}\psi_\varepsilon(\partial_i u)D\varphi\left(\abs{D u}^{p-2}\abs{D\partial_i u}^2+(p-2)\abs{D u}^{p-4}\langle D\partial_i u,D u\rangle^2\right)\,\bd x\\
        =& I+II.
    \end{align*}
    In order to estimate the quanty above, we begin with $I$, which can be estimate using \eqref{Eq:lemmabvinequality}, to get
    \[
        \abs{I}\le C\int_{B_r}\abs{D}\varphi\, \bd x.
    \]
    For $II$,  Lemma \ref{lemma_l2estimate_hessian} yields
    \[
        \abs{II}\le (p-1)\int_{B_r}\abs{D\varphi}\abs{D u}^{p-2}\abs{D^2u}\,\bd x\le C\int_{B_r}\abs{D\varphi}\,\bd x.
    \]
    And by combining this two estimates, the proof is complete.
\end{proof}

\begin{remark}\label{Remark}
Notice that the fact of the inequality \eqref{Eq:lemmabvinequality} has to be understood in the BV-sense it is equivalent to the fact of $\Gamma^{\pm}(u)$ is a set of finite perimeter. Moreover, we known that the sets $\Gamma^{\pm}(u)\cap\{\abs{D u}>0\}$ is locally a $\mathcal{C}^1$-manifolds, hence
\[
\mathcal{H}^{n-1}\left(\Gamma^{\pm}(u)\cap\{\abs{D u>0}\cap B_{1/2}\}\right)\le C,
\]
for some universal constant $C$. Thus we now focus in the set $\Gamma(u)$, the part of the free boundary where the gradient vanishes.
\end{remark}

\begin{lemma}\label{lemma:measureestimate}
Assume that the hypothesis of Theorem \ref{Thm:FinitePerimeter} are in force, and take $r_o$ as in Corollary \eqref{Cor:growth_gradient} . Then, given $\varepsilon_o>0$ we have for any $r<r_o$ 
\[
\frac 1\varepsilon\abs{\{0<\abs{D u}<\varepsilon^{\frac 1{p-\gamma}}\}}\le C.
\] 
Here $C$ is a universal constant.
\end{lemma}
\begin{proof}
	For given $\varepsilon_o>0$ take $r<r_o$. Then for a non-negative $\varphi\in\mathcal{C}_0^{\infty}(B_r)$ such that $\varphi\equiv 1$ on $B_{r/2}$,
	\begin{align*}
		\frac 1{\varepsilon}\abs{\{0<D u<\varepsilon_o^{\frac 1{p-\gamma}}\}}
		\le& C\frac 1\varepsilon\sum_\ell\int_{\{0<\abs{\partial_\ell u}<\varepsilon^{\frac 1{p-\gamma}}\}\cap B_{r/2}}(\Delta_p u)^2\,\bd x\\
		\le & C\frac 1\varepsilon\int_{\{0<\abs{\partial_\ell u}<\varepsilon^{\frac 1{p-\gamma}}\}\cap B_{r/2}}\varphi\left(\abs{D u}^{p-2}\abs{\Delta u}+(p-2)\abs{D u}^{p-4}\abs{u_{x_i}u_{x_j}u_{x_ix_j}}\right)^2\, \bd x\\
		\le& C\int_{\{0<\abs{\partial_\ell u}<\varepsilon^{\frac 1{p-\gamma}}\}\cap B_{r/2}}\varphi\psi_\varepsilon'(\partial_\ell)\abs{D u}^{p-2}\abs{D\partial_\ell u}^2\,\bd x\\
		\le& C\int_{B_r}\abs{D \varphi}\, \bd x.
	\end{align*}
	where in the last two inequalities we have use the fact that
	\[
	\psi_\varepsilon'(t)=\frac{\abs{t}^{p-\gamma}}{\varepsilon}\chi_{\{\abs{t}<\varepsilon^{\frac1{p-\gamma}}\}},
	\]
	and the Lemma \ref{lemabvinequality}.
\end{proof}


The proof of this theorem is a "routinely" application of the  Lemma \ref{lemma:measureestimate}, Corollary \ref{Cor:growth_gradient} and Propositions and \ref{Prop:nondegeneracygamma01}, which we include here for sake of completeness. 
\begin{proof}[Proof of Theorem \ref{Thm:FinitePerimeter}]
 By covering arguments and Remark \ref{Remark:homogenityclass} it is sufficient to prove the theorem for $\Gamma(u) \cap B_{r_0}$ where $r_0$ are as in Lemma \ref{lemma:measureestimate}.

%
Because of Proposition \ref{Prop:nondegeneracygamma01} and Corollary \ref{Cor:growth_gradient} there is a $C$, universal constant, such that for all $\varepsilon>0$ we have
\[
B_{C \varepsilon} \subset B_{\varepsilon} \cap\left\{0<|D u|<\varepsilon^{\frac 1{p-\gamma}}\right \} .
\]

Then
\[
\mathcal{H}^{n-1}\left(\Gamma(u) \cap B_{r_0}\right) \leq C \liminf _{\varepsilon \rightarrow 0} \sum_i \varepsilon^{n-1} .
\]

Finally, let $\{B^i_\epsilon(x_i)\}_{i\in\Lambda}$ a finite cover  of $B_{r_o} \cap \Gamma(u)$, with  centered at $\Gamma(u)$ with at most $N$ overlaps. Using \ref{Eq:lemmabvinequality} and Lemma \ref{lemma:measureestimate} we obtain for $\varepsilon$ small enough
\[
\frac{1}{\varepsilon} \sum_i \varepsilon^n \leq \frac{1}{\varepsilon} C \sum_i\left|B^i \cap\left\{0<|D u|<\varepsilon^{\frac1{p-\gamma}}\right\}\right| \leq \frac{1}{\varepsilon} C N\left|B_{r_0} \cap\left\{0<|D u|<\varepsilon^{\frac1{p-\gamma}}\right\}\right| \leq C .
\]

Therefore,
\[
\mathcal{H}^{n-1}\left(\Gamma(u) \cap B_{r_0}\right) \leq C,
\]
hence, the proof of the theorem is complete.
\end{proof}

\section{The two-dimensional case}

We separate the two-dimensional case because, in this scenario, we know a quantitative estimate for regularity of $p$-harmonic functions and thus we can apply an Liouville-type result to obtain the optimal growth result for $p> 2$. We begin by defining this quantitative estimate. 

\begin{definition}\label{def:Alphap}
Let $p>2$ and $v\in\mathcal{C}(B_1)$ be a weak solution to
\[
\Delta_p v=0\hspace{.2in}\text{in}\hspace{.1in}B_1.
\]
We define $\alpha_p\in(0,1)$ to be the biggest constant such that: there exists  universal constant $C$ such that
\[
\|v\|_{\mathcal{C}^{1,\alpha_p}(B_{1/2})}\le C\|v\|_{L^{\infty}(B_1)}
\]
\end{definition}
\begin{remark}
Notice that the above definition does not depend on the  dimension, although in the this case Berenstein II and Kovalev in  \cite{KB}, show that
$$
\alpha_p> \frac{1}{2p}\left(-3-\frac{1}{p-1}+\sqrt{33+\frac{30}{p-1}+\frac{1}{(p-1)^2}}\right).
$$
Thus, given $p>2$ there exists $\delta_p>0$ such that, for $0<\gamma< 1+\delta_p$
$$
1+\alpha_p>\frac{p}{p-\gamma}.
$$
In the scenario above, we are able to establish the optimal growth result, but now for any $p>2$.   
\end{remark}
Before to proceed with the proof of Theorem \ref{thm:Optimal_growth_dim2} we will prove the following two technical Lemmas.
\begin{lemma}\label{liouville}
   Let $u\in C^{1,\alpha}(B_1)$  be a solution of 
   $$
   \Delta_pu=0,
   $$
   with $u(0)=|Du(0)|=0$ and 
   $$
   |u|\leq C|x|^{1+\alpha}
   $$
   for $0<\alpha<\alpha_p$, then $u=0$.
\end{lemma}
\begin{proof} For $x\in B_1$ consider
$$
v_l(x)=\frac{u(lx)}{l^{1+\alpha}}.
$$
We will show that $\sup_{B_1}|v_l|\leq C$ for sufficiently large $l$. Let $x_l$ be such that
$$
\sup_{B_1}|v_l|=|v(x_l)|.
$$
We have that for $x\in B_1$
$$
|lx|^{1+\alpha}\leq l^{1+\alpha},
$$
Thus, if $|lx|^{1+\alpha}\rightarrow \infty$ as $l\rightarrow \infty$ , by the boundedness assumption, 
$$
|v(x_l)|=\frac{|u(lx_l)|}{l^{1+\alpha}}\leq \frac{|u(lx_l)|}{|lx_l|^{1+\alpha}}\leq C.
$$
On the other hand, if $|lx|^{1+\alpha}\leq C$ as $l\rightarrow \infty$, since $u$ is continuous,
$$
|v(x_l)|=\frac{|u(lx)|}{l^{1+\alpha}}\leq \frac{|C_1|}{|lx|^{1+\alpha}}\leq C.
$$
Moreover, 
$$
\triangle_pv_l=0
$$
thus, by the regularity theory for the $p-$Laplacian ,
$$
|v_l(x)|\leq C|x|^{1+\alpha_p} \quad \mbox{em} \quad B_{1/2}.
$$
Therefore,
$$
\sup_{B_{l/2}}\frac{|u(x)|}{|x|^{1+\alpha_p}}=\sup_{B_{1/2}}\frac{|u(lx)|}{|lx|^{1+\alpha_p}}=\sup_{B_{1/2}}\frac{|u(lx)|}{l^{1+\alpha_p}|x|^{1+\alpha_p}}.
$$
Multiplying and dividing by $l^{1+\alpha}$ we obtain that
$$
\sup_{B_{l/2}}\frac{|u(x)|}{|x|^{1+\alpha_p}}=l^{\alpha-\alpha_p}\sup_{B_{1/2}}\frac{|v_l(x)|}{|x|^{1+\alpha_p}}
$$
thus,
$$
\sup_{B_{l/2}}\frac{|u(x)|}{|x|^{1+\alpha_p}}\rightarrow 0
$$
as $l\rightarrow \infty$, which implies that $u=0$.
\end{proof}
\begin{lemma}\label{caracterizacao de c1alpha}
Let $u\in C^{1,\alpha}(B_1)$, where $0<\alpha<1$. If 
$$
u(0)=|Du(0)|=0
$$
and
$$
\sup_{0<r<1/2}r^{\alpha-\beta}[u]_{C^{1+\alpha}(B_r)}\le A,
$$
for $\beta>\alpha$ and a constant $A>0$, then
$$
|u(x)|\leq A|x|^{1+\beta}, \quad x\in B_{1/2}.
$$
\end{lemma}
\begin{proof}
As $u(0)=0$, by the mean value theorem, for $x\in B_{1/2}$, one has
\begin{equation}\label{ineq1}
   |u(x)|=|u(x)-u(0)|\le|Du(\xi_1)||x|,
\end{equation}
for some $\xi_1\in B_{|x|}$. On the other hand, since $Du(0)=0$ and $Du\in C^{\alpha}\left(B_{1/2}\right)$, we have
$$
\frac{|Du(\xi_1)||x|}{|x|^{1+\beta}}=\frac{|Du(\xi_1)-Du(0)||x|}{|x|^{1+\beta}}\le|x|^{\alpha-\beta}\left[u\right]_{C^{1+\alpha} \left( B_{|x|}\right)}\le A,
$$
which combined with \eqref{ineq1} concludes the proof.
\end{proof}
We are thus led to follow tho the proof of Theorem \ref{thm:Optimal_growth_dim2}.
\begin{proof}[Proof of Theorem \ref{thm:Optimal_growth_dim2}]
Since $u$ is bounded, we have that $u\in C^{1,\alpha}(B_{1/2})$. Moreover,
$$
\|u\|_{C^{1,\alpha}(B_{1/2})}\le C\|u\|_{L^\infty(B_1)},
$$
for a constant $C>0$ depending only on $p$ and $n$. With no loss of generality, we may assume $y=0$ and we need to show that
\begin{equation}\label{nu=1}
    |u(x)|\le C|x|^\frac{p}{p-\gamma}.
\end{equation}
We argue by contradiction and assume that \eqref{nu=1} fails. Then there exist sequences $u_k\in \mathcal{P}_1(M,p,\gamma,y)$ and points $x_k\in B_1$ with 
\begin{equation}\label{branching point uk}
    u_k(0)=|Du_k(0)|=0,
\end{equation}
$$
[u_k]_{C^{1+\alpha}(B_{1/2})}\le2C,
$$
such that
\begin{equation}\label{contradictory assumption}
    |u_k(x_k)|>k|x_k|^{\frac{p}{p-\gamma}}.
\end{equation}
Set 
\begin{equation*}\label{definitiontheta}
    \theta_k(r^{\prime} ):=\sup_{r^{\prime} <r<1/2}r^{1+\alpha-\frac{p}{p-\gamma}}\left[u_k\right]_{C^{1+\alpha}(B_r)}.
\end{equation*}
Observe that $1+\alpha-\frac{p}{p-\gamma}<0$ and
$$
\lim_{r^{\prime}\rightarrow 0}\theta_k(r^{\prime} )=\sup_{0<r<1/2}r^{1+\alpha-\frac{p}{p-\gamma}}\left[u_k\right]_{C^{1+\alpha}(B_r)}.
$$
Then \eqref{contradictory assumption} and Lemma \ref{caracterizacao de c1alpha} yield
$$
\lim_{r^{\prime}\rightarrow 0}\theta_k(r^{\prime})>k.
$$
Thus, there exists a sequence $r_k>\frac{1}{k}$, such that
\begin{equation}\label{contradicao2}
r_k^{1+\alpha-\frac{p}{p-\gamma}}\left[u_k\right]_{C^{1+\alpha}(B_{r_k})}\ge\frac{1}{2}\theta_k(1/k)\ge\frac{1}{2}\theta_k(r_k). 
\end{equation}
From the first inequality in \eqref{contradicao2}, we conclude that 
$$
r_k\rightarrow 0,\,\,\textrm{ as }\,\,k\to\infty.
$$
Set now 
\begin{equation}\label{definitionvk}
    v_k(x):=\frac{u_k(r_k x)}{\theta_k(r_k)r_k^{\frac{p}{p-\gamma}}}
\end{equation}
and note that, for $1\le R\le\frac{1}{2r_k}$, one has
\begin{equation}\label{vkestimate}
\begin{split}
    \left[v_k\right]_{C^{1+\alpha}(B_R)}&=\frac{1}{\theta_k(r_k)r_k^{\frac{p}{p-\gamma}}}\left[u_k\right]_{C^{1+\alpha}(B_{r_kR})}r_k^{1+\alpha}\\
    &=\frac{(r_kR)^{1+\alpha-\frac{p}{p-\gamma}}}{\theta_k(r_k)}\left[u_k\right]_{C^{1+\alpha}(B_{r_k R})}R^{\frac{p}{p-\gamma}-1-\alpha}.
    \end{split}
\end{equation}
Since
$$
(r_k R)^{1+\alpha-\frac{p}{p-\gamma}}\left[u_k\right]_{C^{1+\alpha}(B_{r_k R})}\le\theta_k(r_k R)\le\theta_k(r_k),
$$
employing \eqref{vkestimate}, we get
\begin{equation}\label{seminorm estimate on vk}
    \left[v_k\right]_{C^{1+\alpha}(B_R)}\le R^{\frac{p}{p-\gamma}-1-\alpha},\quad\forall R\ge1.
\end{equation}
Now, if $\eta$ is a smooth function that is $1$ in $B_{1/2}$ and $0$ outside $B_1$, then for any unit vector $e$, one has
$$
\int_{B_1}\eta\cdot D_e v_k\,dx=\int_{B_1}D_e\eta\cdot v_k\,dx\leq C_n,
$$
where $C_n>0$ is a constant depending only on $n$. Therefore, there exists $z\in B_1$ such that $|Dv_k(z)|\leq C_n$, and \eqref{seminorm estimate on vk} implies 
$$
|D_ev_k(x)-D_ev_k(z)|\le R^{\frac{p}{p-\gamma}-1-\alpha}|x-z|^{\alpha}.
$$
Hence, for $x\in B_R$ and $1\leq R \le\frac{1}{2r_k}$, we get
$$
|D_ev_k(x)|\le C_n+ R^{\frac{p}{p-\gamma}-1-\alpha}|x-z|^{\frac{p}{p-\gamma}}\le CR^{\frac{p}{p-\gamma}-1},
$$
for a constant $C>0$ depending only on $n$. Moreover
$$
|v_k(x)|=|v_k(x)-v_k(0)|\leq |Dv_k(\xi)||x|\leq C|x|^{\frac{p}{p-\gamma}-1}|x|=C|x|^{\frac{p}{p-\gamma}}.
$$
Thus, up to a subsequence, $v_k$ converges to some $v_0$ in $C^{1+\alpha}_{\text{loc}}(B_R)$, as $k\to\infty$, and thanks to the last inequality,
\begin{equation*}\label{v0estimate}
\left[v_0\right]_{C^{1+\alpha}(B_R)}\le R^{\frac{p}{p-\gamma}-1-\alpha}, \quad \forall R\ge1.
\end{equation*}
Note that \eqref{branching point uk} and \eqref{definitionvk} imply
\begin{equation*}
    v_0(0)=|D v_0(0)|=0,
\end{equation*}
therefore
$$
|v_0(x)|\leq C|x|^{\frac{p}{p-\gamma}}.
$$
Additionally, from the second inequality in \eqref{contradicao2}, we deduce that
\begin{equation}\label{limitseminormestimate2}
     \left[v_0\right]_{C^{1+\alpha}(B_{1})}\ge\frac{1}{2}.
\end{equation}

\begin{equation*}\label{becka}
    \Delta_p v_k(x)
    =\frac{1}{(\theta_k(r_k))^{p-1}r_k^{\frac{p(\gamma-1)}{p-\gamma}}}\left[\left(u_k(x)\right)_+^{\gamma-1}-\left(u_k(x)\right)_-^{\gamma-1}\right]=\frac{1}{(\theta_k(r_k))^{p-\gamma}}\left[\left(v_k(x)\right)_+^{\gamma-1}-\left(v_k(x)\right)_-^{\gamma-1}\right].
\end{equation*}
We can then pass to the limit, as $k\to\infty$, and arrive at
\begin{equation*}\label{limit equation}
    \Delta_pv_0(x)=0\,\,\mbox{ in }\,\,\mathbb{R}^n.
\end{equation*}
Note that \eqref{branching point uk} and \eqref{definitionvk} imply
\begin{equation*}\label{branching point v0}
    v_0(0)=|D v_0(0)|=0,
\end{equation*}
thus, by Lemma  \ref{liouville}, we have $v_0=0$ which contradicts \eqref{limitseminormestimate2}, thereby proving the result.

\end{proof}
\begin{remark} We which to conclude this section with two important remarks. The first one is that in two dimension the proof of Theorem \ref{Thm:FinitePerimeter} remains valid, since $p>2$, thus we can claim that $\Gamma(u)$ is a finite perimeter set. The second remark is to highlight the fact that the crucial step is the effective quantification on $\alpha_p$ in two dimension, thus if known the behavior of $\alpha_p$ in greatest dimensions, our reasoning still valid.
\end{remark}

\bigskip

\noindent{\bf Acknowledgements:} J. Correa  has received partial support from CNPq-Brazil under Grant No. 408169/2023-0. D. dos Prazeres has received partial support from CNPq-Brazil under Grant No. 305680/2022-6.\\



\bigskip
	\nocite{*}
	\bibliographystyle{plain}
	\bibliography{bibliography}

\bigskip

\noindent\textsc{Julio C. Correa-Hoyos}\\
Instituto de Matem\'atica e Estat\'istica\\Universidade do Estado do Rio de Janeiro\\20550-013, Maracan\~a, Rio de Janeiro - RJ, Brazil.\\\noindent\texttt{julio.correa@ime.uerj.br}.
\bigskip

\noindent\textsc{Disson dos Prazeres }\\
 Departamento de Matemática\\ 
 Universidade Federal de Sergipe\\
\noindent\texttt{disson@mat.ufs.br}

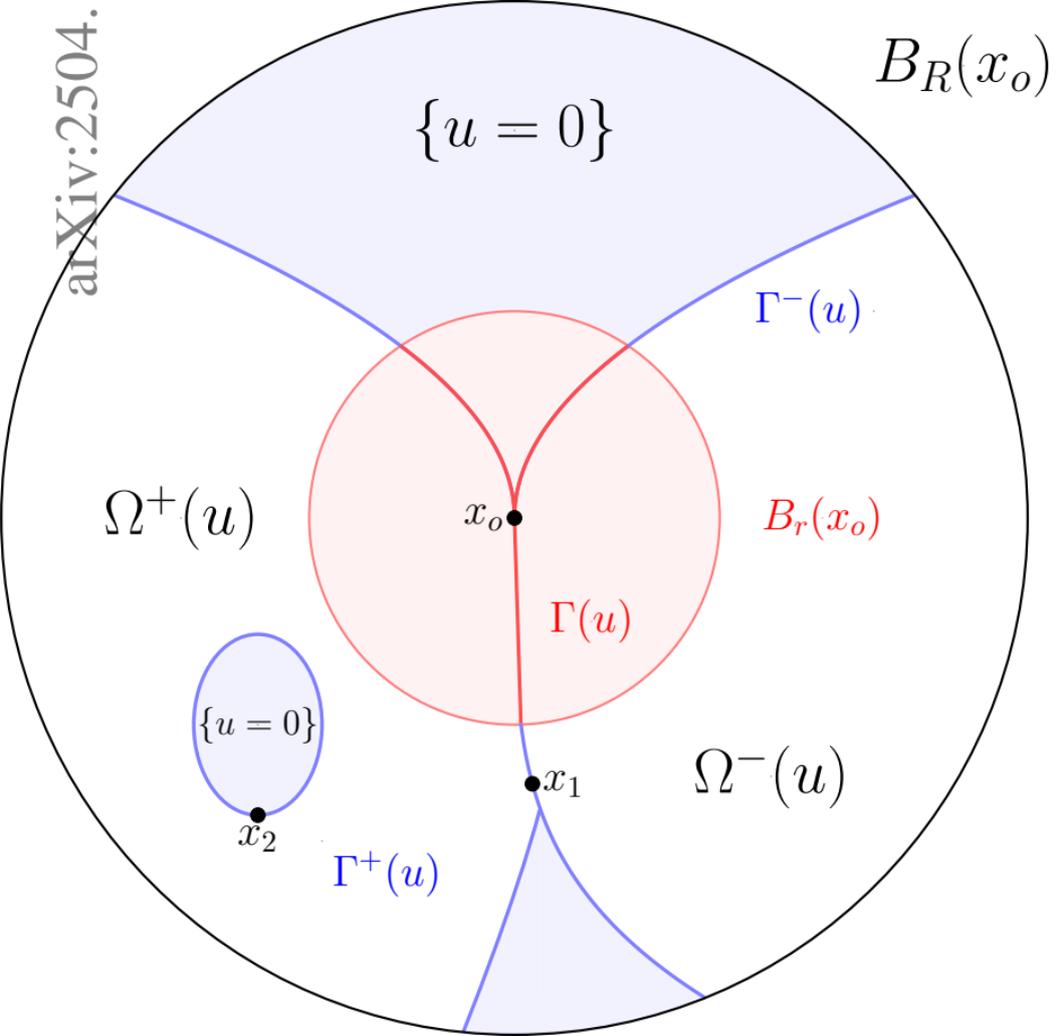
\end{document}